\documentclass[12pt]{amsart}

\usepackage{amsmath,amssymb,amsthm}

\numberwithin{equation}{section}
\theoremstyle{plain}
\newtheorem{theorem}[equation]{Theorem}

\newtheorem{prop}[equation]{Proposition}

\newtheorem{lem}[equation]{Lemma}
\newtheorem{lemma}[equation]{Lemma}

\newtheorem{corollary}[equation]{Corollary}

\theoremstyle{remark}
\newtheorem{remark}[equation]{Remark}

\theoremstyle{definition}
\newtheorem{definition}[equation]{Definition}

\newtheorem{question}[equation]{Question}

\newcommand{\lra}{\longrightarrow}
\newcommand{\ra}{\rightarrow}
\newcommand{\restr}{\mbox{\Large \(|\)\normalsize}}
\newcommand{\N}{\mathbb N}
\newcommand{\R}{\mathbb R}

\newcommand{\E}{\mathbb E}

\renewcommand{\P}{\mathbb P}

\newcommand{\acts}{\curvearrowright}
\newcommand{\aut}{\operatorname{Aut}}
\renewcommand{\bot}{\operatorname{bot}}
\newcommand{\ch}{\operatorname{CH}}

\newcommand{\fix}{\operatorname{Fix}}

\newcommand{\isom}{\operatorname{Isom}}
\newcommand{\rank}{\operatorname{Rank}}
\newcommand{\sh}{\operatorname{Sh}}
\newcommand{\slope}{\operatorname{slope}}

\renewcommand{\top}{\operatorname{top}}
\newcommand{\tits}{\partial_T}

\def\D{\partial}
\newcommand{\al}{\alpha}

\def\De{\Delta}

\def\eps{\epsilon}
\def\ga{\gamma}

\def\Ga{\Gamma}

\def\La{\Lambda}
\def\lang{\langle}
\def\<{\lang}
\def\>{\rangle}
\def\lra{\longrightarrow}

\def\mod{\De_{mod}}

\def\ol{\overline}

\def\ra{\rightarrow}

\def\si{\sigma}
\def\Si{\Sigma}

\def\geo{\partial_{\infty}}

\def\tangle{\angle_{T}}
\def\defeq{:=}

\def\BI{\begin{itemize}}
\def\EI{\end{itemize}}

\begin{document}

\title{Rigidity of invariant convex sets in symmetric spaces}
\author{Bruce Kleiner}
\thanks{Supported by  NSF grant DMS-0204506.}
\author{Bernhard Leeb}
\date{\today}
\maketitle

\begin{abstract}
The main result implies
that a proper convex subset of an irreducible higher rank
symmetric space cannot have Zariski dense stabilizer.
\end{abstract}

\section{Introduction}
In this paper we study convex subsets of symmetric spaces,
and their stabilizers.  The main results show that in the
higher rank case convex sets are strongly restricted, and 
under mild assumptions can only arise from rank $1$ constructions.
This rigidity phenomenon for convex
subsets is yet another example of a rigidity
property enjoyed by higher rank symmetric spaces that
has no analog for rank $1$ symmetric spaces.

One can generate a
supply of convex subsets of any Hadamard space
by starting with geodesic segments, geodesic rays, complete
geodesics,  and horoballs, and then taking tubular neighborhoods and
intersections.
When $X$ is a Hadamard manifold with pinched negative
curvature convex subsets are abundant: 
by a theorem of Anderson \cite{anderson},
any closed subset $A$ of the geometric
boundary $\geo X$ is the limit set of a closed convex
subset $Y\subset X$.  On the other hand, for general Hadamard
spaces (or manifolds) it can be difficult to control
the convex hull of even ``small'' subsets, like
the union of three rays.

A  group $\Ga$ of isometries of a Hadamard space $X$
is {\em convex cocompact} if there is a $\Ga$-invariant
convex subset $C\subset X$ with compact quotient $C/\Ga$.
Discrete convex cocompact  subgroups of the isometry group of 
hyperbolic $3$-space are an important class in the theory 
of Kleinian groups; basic examples are uniform lattices,
Schottky groups
and quasi-Fuchsian groups.  Analogous examples exist
in $\isom(H^n)$, as well as the isometry groups of 
other rank $1$ symmetric spaces.   In a higher
rank symmetric space of noncompact type, 
one can produce examples by taking products of
 uniform lattices and  
rank $1$ convex cocompact groups.  In 1994,  Corlette
asked if this was essentially the only way to produce
discrete convex cocompact groups.  The answer  is 
 yes, see Theorem \ref{mainconvexcocompact} below; in fact
the theorem is proved by reducing it to the case of convex
subsets with Zariski dense stabilizer:

\begin{theorem}
\label{mainzdense}
Let $X=\E^n\times Y$, where $Y$ is  a symmetric space 
of noncompact type, and
let $X=\E^n\times Y_1\times Y_{\geq 2} $ denote the decomposition
of $X$ into the Euclidean factor, the
product of the irreducible rank $1$ factors, the product of the
higher rank factors.
Suppose $\Ga\subset  \isom(X)=\isom(\E^n)\times\isom(Y)$
is a subgroup
whose  projection to $\isom(Y)$
is Zariski dense in the identity component $\isom_o(Y)$,
and whose projection to
$\isom(\E^n)$ does not preserve a proper affine subspace of $\E^n$.
If $C\subset X\defeq \E^n\times Y$
is a $\Ga$-invariant closed convex set,
then  $C=\E^n\times C_1\times Y_{\geq 2}$,
where $C_1\subset Y_1$ is a closed convex subset.
Furthermore,  for each de Rham factor $X_i$ of
$Y_1$, there
is a  $\Ga$-invariant subset $\hat C_i\subset X_i$  such that
\begin{itemize}
\item $\hat C_i$ is the 
closed convex hull of its limit set.

\item $|\geo C_i|=\infty$,

\item $\hat C_1\defeq \prod_i\hat C_i\subset C_1$.

\item  $\geo\hat C_1=\geo C_1$. 
\end{itemize}
\end{theorem}
We recall that by convention, a symmetric space of noncompact
type has no Euclidean de Rham factor.  Note that a subgroup of
$\isom_o(Y)$ is Zariski dense if and only if it neither
fixes a point in the Tits  boundary $\tits Y$  nor preserves
a proper symmetric subspace of $Y$.

\begin{corollary}
If $X$ is a symmetric space of noncompact type
with no rank $1$ de Rham factors and $\Ga\subset\isom_o(X)$
is a Zariski dense subgroup, then $X$ contains no
proper closed $\Ga$-invariant convex subsets.
\end{corollary}

For discrete convex cocompact groups, we have the following
structural result:

\begin{theorem}
\label{mainconvexcocompact}
Let $X=\E^n\times Y$, where $Y$ is a symmetric space of 
noncompact type.  Suppose
$\Ga\subset \isom(X)= \isom(\E^n)\times\isom(Y)$ 
is a discrete subgroup acting cocompactly on a closed convex subset
$C\subset X$, and assume $\Ga$ does not preserve any proper symmetric
of $X$.   Then 
$\Ga$ projects to a subgroup of 
$\isom(Y)$ which is Zariski dense in $\isom_o(Y)$, 
and the conclusions of Theorem \ref{mainzdense}
apply to $C$.
\end{theorem}
If  a convex cocompact subgroup
$\Ga\subset \isom(X)$ preserves a proper symmetric subspace
$Z\subset X$, then it acts convex cocompactly
on $Z$ -- just intersect a sufficiently big tubular neighborhood
of a $\Ga$-invariant convex set with $Z$.  Therefore there is 
no loss of generality in assuming $X$ contains no proper
$\Ga$-invariant symmetric subspace.

\begin{corollary}
If $X$ is a symmetric space of noncompact type with no
rank $1$ de Rham factors, and $\Ga\subset\isom(X)$
is a discrete subgroup acting cocompactly on a closed
convex subset $C\subset X$, then either $C=X$ and 
$\Ga$ is a uniform lattice in $\isom(X)$, or
$\Ga$ preserves a proper symmetric subspace of $X$.
\end{corollary}

We give a brief outline of the proof of Theorem \ref{mainzdense}
in the case the Euclidean factor is absent, and $Y$
is an irreducible higher rank symmetric space.  The first
step is to apply a Theorem of Benoist \cite{benoist}, which
implies one may find an open neighborhood $U$ of a pair
of antipodal points $\xi_1,\xi_2$ in the  Tits boundary
$\tits X$, such that $U$ is contained in the limit set
of $\Ga$.  Applying a result from \cite{klle}, we deduce
that the geometric boundary  of $C$ is a top dimensional subbuilding
$B$ of the Tits boundary of $X$, which  is a closed
subset with respect to the topology of the geometric boundary 
$\geo X$.  The main 
step in the paper, implemented in Theorem \ref{filling},
 is to show that any such building is contained
in the geometric boundary of
a proper symmetric subspace $Y$, unless it coincides with
$\tits X$;  the Zariski density assumption
rules out the former possibility in the case at hand.  We remark that Theorem
\ref{filling} applies to products of symmetric spaces and
Euclidean buildings, and may be of independent interest.

In view of the  results in this paper
one may wonder 
whether sufficiently large convex sets in symmetric spaces 
of noncompact type 
or in spherical buildings (such as Tits boundaries of symmetric spaces) 
are rigid. 

\begin{question}
Suppose $C\subset B$ is a convex subset of a spherical building.
If $C$ does not have circumradius $\leq \frac{\pi}{2}$,
must $C$ itself be a spherical building? 
\end{question}

It is unclear what one should expect here. 
A.\ Balser and A.\ Lytchak \cite{ballyt}
proved a partial result regarding convex subsets 
invariant under a group action,
namely if $\dim(C)\leq2$ and $C$ is not a spherical building
then $\isom(C)$ has a fixed point in $C$. 

After the first version of this paper was written,
Quint informed the authors of very interesting
related work \cite{quint} on Zariski dense subgroups of 
semi-simple groups.  His paper addresses an alternate
definition of convex cocompact groups which is equivalent
to the usual definition for rank $1$ symmetric spaces
but differs from ours in the higher rank case; 
for this reason  it is difficult to make a direct
comparison between the results of \cite{quint}
and the theorems above.  We mention that his main 
result also applies to discrete subgroups of semi-simple $p$-adic
groups.

 The authors proved slightly weaker
versions of Theorems \ref{mainzdense} and \ref{mainconvexcocompact}
in 1998, and spoke publicly about them in the subsequent year.

\tableofcontents

\section{Preliminaries}

\subsection{Hadamard spaces}
We recommend \cite{ball,brha,klle} as  references for
Hadamard space facts.  

The term {\em Hadamard space} is a synonym for a $CAT(0)$-space.  

If $X$ is a Hadamard space, we denote the geometric 
boundary by $\geo X$, the Tits boundary by $\tits X$, 
and the Tits angle between $\xi_1,\xi_2\in\tits X$
by $\tangle(\xi_1,\xi_2)$.

Recall that the set underlying $\geo X$ is the set of asymptote
classes of geodesic rays, and that this  may be identified with
the set of rays leaving a given basepoint $p\in X$.
If $x_1,x_2\in X$, $Y\subset X$ is a subset, $y_i\in Y$ is a sequence
with $\lim_{i\ra\infty}d(y_i,p)=\infty$,
then the segments $\ol{x_1y_i}$ converge to a ray
$\ol{x_1\xi}$ iff the segments $\ol{x_2y_i}$
converge to a ray $\ol{x_2\xi}$.  Thus the set of
rays which can be obtained as limits in this fashion,
as $\{y_i\}$ ranges over all such sequences, is a collection
of asymptote classes and therefore determined a subset
of $\geo X$, the {\em limit set of $Y$}, which we denote
 by $\La(Y)$.  

\begin{lemma}
\label{lemconvexcone}
If $C\subset X$ is a closed convex subset, and 
$p\in C$, then every ray
$\ol{p\xi}$ is contained in $C$, for $\xi\in \La(C)$.
\end{lemma}
\proof 
This follows from the convexity of $C$ and
the definition of the limit set,
since we are at liberty select the basepoint.
\qed

\begin{definition}
\label{defconvex}
A subset $Y$ of a $CAT(1)$ space $Z$ is {\em convex}
if it contains every segment of the form
$\ol{\xi_1\xi_2}$, where $\xi_1,\xi_2\in Y$
and $d_Z(\xi_1,\xi_2)<\pi$.
\end{definition}

\begin{lemma}
\label{lemtitsconvex}
Let $X$ be a proper Hadamard space, and let $C\subset X$
be a closed convex subset.  Then the limit set of $C$
in $\geo X$ determines a convex subset of $\tits X$,
which is isometric to the Tits boundary of $C$, viewed
as a Hadamard space.
\end{lemma}
\proof
The isometric embedding
$C\ra X$ of Hadamard spaces 
induces an isometric embedding $\tits C\ra \tits X$
of Tits boundaries.  Since $\tits C$ is a $CAT(1)$ space,
the image of the embedding is convex. 
\qed

\begin{lemma}
\label{prelimrsplit}
If $\Ga\acts X$ is a discrete, cocompact, isometric
action on a Hadamard space $X$,  
 and $\Ga$ fixes a point
$\xi\in\tits X$, then there is a geodesic $\ga\subset X$
such that $\xi\in\tits\ga$ and the parallel set
$\P(\ga)\subset X$ is $\Ga$-invariant.
\end{lemma}
\proof
We may assume that $X$ contains no proper, closed, convex,
$\Ga$-invariant nonempty subset, by applying Zorn's lemma.

Note that any element $g\in Z(\Ga)$ is semi-simple
and its minimum displacement set, $\min(g)\subset X$,
 is a closed, convex, and $\Ga$-invariant
subset; therefore by assumption we have $\min(g)=X$.
Thus elliptic elements in $Z(\Ga)$ act trivially on 
$X$ and nonelliptic elements act by Clifford translations,
i.e. they translate along the $\R$-factor of a product 
splitting $X=\R\times Z$.  Hence $X$ admits a product
structure
\begin{equation}
\label{eqnclifford}
X=\E^n\times Y
\end{equation}
where $Z(\Ga)$ acts by translations on $\E^n$ and trivially
on $Y$.

Pick $p\in X$, and a finite generating set $\Si\subset \Ga$.
Let $C\defeq \max_{\si\in\Si}d(\si p,p)$.  Note that the
ray $\ol{p\xi}\subset X$ lies in the closed convex set
$$
\De\defeq\{x\in X\mid \mbox{For all $\si\in \Si$, $d(\si x,x)\leq C$}\},
$$
since for all $g\in \Ga$ and every $x\in \ol{p\xi}$,
we have $d(gx,x)\leq d(gp,p)$ because $\ol{p\xi}$ and 
$\ol{(gp)\xi}$ are asymptotic rays.
By a standard argument the centralizer, $Z(\Si)=Z(G)$, of the set
$\Si$ acts cocompactly on $\De$, which implies that
$\ol{p\xi}$ is contained in a finite tubular neighborhood
of an $n$-flat $\E^n\times\{y\}$ of the product decomposition
(\ref{eqnclifford}).  Hence $\xi\in \tits \E^n$, and this 
implies the lemma.
\qed

\subsection{Affine and  concave functions on convex sets}

\begin{lemma}
\label{concavecomplete}
Let $Z$ be a geodesic metric space with extendible geodesics.
Then any concave function $Z\to[0,\infty)$ is constant.
\end{lemma}
\proof Trivial. \qed

\begin{lemma}
\label{nonnegconcrank1}
Let $Z$ be a CAT(-1)-space whose ideal boundary 
$\geo Z$ consists of at least two points. 
Suppose that there is no proper closed convex subset of $Z$ 
whose ideal boundary equals $\geo Z$.
Then any continuous concave function $f:Z\to[0,\infty)$ is constant. 
\end{lemma}
\proof
We first observe that $f$ is constant along each complete geodesic.
Furthermore,
$f$ is non-decreasing along each geodesic ray,
and the restriction of $f$ to a compact geodesic segment 
assumes its minimum at one of the endpoints. 

Note that, by assumption,
$Z$ contains at least one complete geodesic.
Let $l$ be a complete geodesic and $z\in Z$ be an arbitrary point.
Denote by $\rho_1,\rho_2:[0,\infty)\to Z$
the rays emanating from $z$ and asymptotic to the two ends of $l$.
Then $f$ is $\geq f(z)$ along each segment connecting 
$\rho_1(t)$ to $\rho_2(t)$ for $t\geq0$. 
Since $Z$ is CAT(-1) these segments converge to the line $l$.
The continuity of $f$ then implies that 
$f(l)\geq f(z)$. 
Thus $f$ assumes on $l$ its maximum
which we denot by $m$. 

It follows that $f$ equals $m$ 
on the union $H_1$ of all lines in $Z$.
Consider the ascending sequence of subsets $H_n\subset Z$
defined inductively by requiring that $H_{n+1}$
is the union of all segments with endpoints in $H_n$. 
Then the sequence of suprema $\sup(f|_{H_n})$ is non-decreasing.
Hence $m=\sup(f|_{H_1})\leq\sup(f|_{H_n})\leq m$
and $f\equiv m$ on the closure of $\cup_{n\in\N} H_n$.
This closure is a closed convex subset of $Z$ with the same ideal boundary
and, by assumption, equals $Z$.
\qed

By an {\em affine} function on a geodesic metric space  
we mean a function whose
restriction to each geodesic segment is an affine function.
\begin{lem}
\label{affinerank1}
Let $Z$ be a CAT(-1)-space whose ideal boundary 
$\geo Z$ consists of at least three points. 
Then any affine  continuous function $f:Z\to\R$ is constant. 
\end{lem}
\proof
We first observe that the slope of $f$ 
along a geodesic ray 
depends only on the ideal point represented by it. 
Indeed,
let $\rho_1,\rho_2:[0.\infty)\to Z$
be two rays parametrized by unit speed. 
Since the geodesic segments connecting $\rho_1(0)$ with $\rho_2(t)$ 
converge to the ray $\rho_1$
it follows using continuity
that the slope of $f$ along $\rho_1$ equals its slope along $\rho_2$. 

Since any two ideal points in $\geo Z$ may be connected 
by a complete geodesic in in $Z$
it follows that the slopes of $f$ at any two ideal points 
have opposite sign.
Since $\geo Z$ contains at least three points
the slopes of $f$ must be zero at all ideal points,
i.e.\ $f$ is constant along every geodesic ray.

The same reasoning as in the proof of Lemma \ref{nonnegconcrank1}
above shows 
that for any point $z$ and any complete geodesic $l$ in $Z$ 
we have $f(z)=f(l)$.
Thus $f$ is constant.
\qed

\begin{lem}
\label{affinehigherrank}
Let $Z$ be a symmetric space of noncompact type and higher rank
without Euclidean de Rham factor.
Then any affine  continuous function $f:Z\to\R$ is constant. 
\end{lem}
\proof
We may apply Lemma \ref{affinerank1} 
to (nonflat) totally geodesic subspaces of rank one
and get that $f$ is constant on any such subspace.

Let $F$ be a maximal flat.
Then $f|_F$ is affine . 
The previous remark implies 
that the gradient of $f|_F$ at a point $z\in F$ 
must be tangent to every singular hyperplane $H$ through $z$
because the lines in $F$ perpendicular to $H$ 
lie in a rank one subspace. 
Since $Z$ has no Euclidean factor
the intersection of all these hyperplanes $H$ is just the point $z$. 
We conclude that $f$ is constant along every maximal flat;
since any two points
lie in a maximal flat, this implies that $f$ is constant
on $Z$. 
\qed

\subsection{Asymptotic slopes of convex functions}
\label{prelimslopes}

Let $Z$ be a Hadamard space 
and $f:Z\to\R$ a continuous convex function.
For a unit speed goedesic ray $\rho:[0,\infty)\to Z$
we define the {\em asymptotic slope} of $f$ along $\rho$ 
as $\slope_f(\rho):=\lim_{t\to\infty}\frac{f(\rho(t))}{t}
\in\R\cup\{\infty\}$. 

\begin{lem}
For any two asymptotic unit speed rays 
$\rho_1$ and $\rho_2$, 
$\slope_f(\rho_1) = \slope_f(\rho_2)$.
\end{lem}
\proof
Since the segments connecting $\rho_2(0)$ with $\rho_1(t)$
Hausdorff converge to $\rho_2$
one estimates using the continuity of $f$ that 
$f(\rho_2(t))\leq C+\slope_f(\rho_1)\cdot t$ 
for $t\geq0$ 
and hence $\slope_f(\rho_2) \leq \slope_f(\rho_1)$.
Symmetry implies equality.
\qed

Thus we may  speak of the asymptotic slope, $\slope_f(\xi)$, 
at an ideal point $\xi\in\geo Z$. 

\begin{lem}
$\slope_f:\geo Z\to \R\cup\{\infty\}$ 
is lower semicontinuous 
with respect to the cone topology.
\end{lem}
\proof
Consider a sequence of unit speed rays $\rho_n$
with same initial point 
which Hausdorff converges to the ray $\rho$.
Since $\slope_f(\rho_n)\geq\frac{f(\rho_n(t))-f(\rho_n(0))}{t}$
for $t\geq0$ 
we obtain 
\[
\liminf_{n\to\infty}\slope_f(\rho_n)\geq\frac{f(\rho(t))-f(\rho(0))}{t}
\buildrel{t\nearrow\infty}\over\lra \slope_f(\rho) .
\]
\qed

As a consequence,
$\slope_f$ attains a minimum if $Z$ is locally compact.

\begin{prop}
\label{propuniqueminimum}
If $\slope_f:\geo Z\to \R\cup\{\infty\}$ 
assumes negative values 
then it has a unique minimum. 
\end{prop}
\proof
Let $\xi_1,\xi_2\in\geo Z$ be ideal points with 
$\slope_f(\xi_i)\leq-a<0$ 
and $\tangle(\xi_1,\xi_2)\geq\eps>0$. 
Let $\rho_i$ be unit speed rays 
emanating from the same point $o\in Z$
and asymptotic to the ideal points $\xi_i$. 
For the midpoints $m(t)$ of the segments 
$\ol{\rho_1(t)\rho_2(t)}$ 
holds 
\[
\limsup_{t\to\infty}\frac{d(o,m(t))}{t}
\leq\cos\frac{\tangle(\xi_1,\xi_2)}{2}
\leq \cos\frac{\eps}{2}
.\]
Moreover $f(m(t))\leq const-at$.
The segments $\ol{om(t)}$ Hausdorff converge to a ray $\mu$
which therefore satisfies
$\slope_f(\mu)\leq-a(\cos\frac{\eps}{2})^{-1}$.

It follows that any sequence $(\xi_n)$ in $\geo Z$ 
with $\slope_f(\xi_n)\searrow\inf \slope_f$
is a Cauchy sequence with respect to the Tits metric.
Hence $\slope_f$ has a unique minimum on $\geo Z$.
\qed

\subsection{Spherical buildings}
We refer the reader to \cite{klle,ronan,tits} for further discussion of the
material here.  

We will be using the geometric definition of 
spherical buildings from \cite{klle}, which we now recall.

Let $(S,W)$ be a spherical Coxeter complex, so $S$ is a Euclidean
sphere and $W$ is a finite group generated by reflections acting
on $S$.  A {\em
spherical building modelled on $(S,W)$} is a
$CAT(1)$-space $B$ together with a collection ${\mathcal A}$ of
isometric embeddings $\iota:S\ra B$, called {\em charts},
which satisfies properties
SB1-2 described below 
and which is closed under precomposition with isometries in $W$. 
An {\em apartment} in $B$ is the image of a chart $\iota:S\ra B$;
$\iota$ is a chart of the apartment $\iota(S)$.

\smallskip\noindent
{\bf SB1: Plenty of apartments.}
Any two points in $B$ are contained in a common apartment.

\smallskip\noindent
Let $\iota_{A_1}$, $\iota_{A_2}$ be charts for apartments $A_1$,
$A_2$, and let $C=A_1\cap A_2$, $C'=\iota_{A_2}^{-1}(C)\subset S$.
The charts $\iota_{A_i}$ are {\em $W$-compatible} if
$\iota_{A_1}^{-1}\circ\iota_{A_2}\restr_{C'}$ is the restriction
of an isometry in $W$.

\smallskip\noindent
{\bf SB2: Compatible apartments.}
The charts are $W$-compatible.

\subsection{Root groups}
\label{prelimrootgroups}

If $B$ is spherical building, and $a\subset B$ is a root, 
then the {\em root group of $a$} is the collection $U_a$ of 
building automorphisms of $B$ which fix $a$ pointwise, as
well as any chamber $\si\subset B$ such that $\si\cap a$
is a panel $\pi$ which is not contained in the wall
$\D a$.  The building $B$ is {\em Moufang} if for every
root $a\subset B$, the group $U_a$ acts transitively on the
set of roots opposite $a$.  

\medskip
{\em Properties of root groups:}

$\bullet$ When all  the join factors of $B$  have 
 dimension at least $1$, then $U_a$ acts freely on the
collection of roots opposite $a$.

$\bullet$ When $X$ is a symmetric space of noncompact
type and $B\defeq\tits X$, then $B$ is a Moufang
building and  $G\defeq\isom_o(X)$
acts effectively on $B$ by building automorphisms,
so we may view $G$ as a subgroup of $\aut(B)$.  Each 
root group of $B$ is contained in $G$, and is a unipotent
subgroup \cite[pp. 77-78]{tits}.
Furthermore, $G$ is generated by the root groups of $B$.

\subsection{Groups acting on symmetric spaces}
\label{prelimsymmspace}

Let $X$ be a symmetric space of noncompact type,
and let $G\defeq\isom_o(X)$.
We will require the following well known facts \cite{mos,boreltits}:
\begin{itemize}
\item
A  subgroup
$H\subset G$ is Zariski dense if and only if
$H$ neither fixes a point in $\tits X$ nor
preserves a proper symmetric subspace.
 
\item A proper subgroup  $H\subsetneq G$ with
finitely many connected components is not
Zariski dense; in particular $H$ must either
fix a point in $\tits X$ or preserve a proper
symmetric subspace.
\end{itemize}

\begin{remark}
If a Zariski dense subgroup of a real simple group
 is not dense in the usual
topology, then it must be discrete.
\end{remark}

\section{Top dimensional subbuildings in the boundary 
of a symmetric space}

In this section we prove:

\begin{theorem}
\label{filling}
Suppose 
\begin{equation}
\label{xdecomp}
X=X_1\times\ldots\times X_k
\end{equation} is a product of
irreducible symmetric spaces of noncompact type,
irreducible Euclidean buildings with
discrete affine Weyl groups, and Euclidean spaces.
Let $B\subset\tits X$ be a top dimensional
subbuilding which is closed with respect to the  topology
of the geometric boundary $\geo X$,
and which is not contained in the boundary 
of any proper subspace  $Y\subset X$ of the form
$Y=Y_1\times\ldots\times Y_k$, where $Y_i\subset X_i$
is either a totally geodesic subspace or a subbuilding,
according to the type of $X_i$.  Then
there is a join decomposition
$$B=B_1\circ\ldots\circ B_k$$
where $B_i\defeq B\cap \tits X_i$, such that
$B_i=\tits X_i$ unless $X_i$ is an irreducible
rank $1$ symmetric space of noncompact type.
\end{theorem}

We begin the proof by observing that if there is more than 
one factor in the product decomposition (\ref{xdecomp}),
then by \cite[Prop. 3.3.1]{klle}, $B$ and $\tits X$ will admit  
corresponding compatible
join decompositions 
$$B=B_1\circ\ldots\circ B_k,\quad 
\tits X=\tits X_1\circ\ldots\circ\tits X_k,$$
and hence it is sufficient to prove the theorem for the irreducible
factors $X_i$ separately.  So henceforth we will assume that
$X$ is irreducible.  If $X$  is Euclidean, then $\tits X$
is the only top dimensional subbuilding of $\tits X$, and 
so this case is trivial.

\subsection{The case when $X$ is a Euclidean building}
Let $Y\subset X$ be the union of the collection of 
apartments $A\subset X$
such that $\tits A\subset B$.   

\begin{lemma}
Any two chambers $\si_1,\si_2\subset Y$ lie in an apartment
$A\subset X$ which is entirely contained in $Y$.
\end{lemma}
\proof
By the definition of $Y$, for $i=1,2$ there exists an apartment
$A_i\subset X$ such that $\tits A_i\subset B$ and 
$\si_i\subset A_i$.  For $i=1,2$, choose an interior point
$p_i\in \si_i$, and consider the geodesic segment $\ol{p_1p_2}\subset X$.
By perturbing $p_2$ slightly, if necessary, we may assume that
the $\mod$-direction of $\ol{p_1p_2}$ is regular.  We may
prolong $\ol{p_1p_2}$ to a complete regular geodesic $\ga\subset X$
by concatenating it with rays $\ol{p_1\xi_1}\subset A_1$,
$\ol{p_2\xi_2}\subset A_2$.  Since $\tits \ga=\{\al_1,\al_2\}$
where $\al_i\subset\tits A_i\subset B$ are regular, there
is a unique apartment $\tits A\subset \tits X$ containing
$\geo \ga$, and it is contained in $B$.  Then by the definition
of $Y$ we have $A\subset Y$, and since $\ga$ is regular
and $\tits \ga\subset \tits A$, we get $\ga\subset A$. 
This implies that $\si_i\subset A$, since $\si_i\cap A\supset \{p_i\}\neq
\emptyset$ is a subcomplex of $X$.
\qed

\medskip
\noindent
The lemma implies that $Y$ is a subbuilding of $X$ with Tits
boundary $B$.  By assumption we must therefore have $B=\tits X$,
which proves Theorem \ref{filling} in this case.

\subsection{$X$ is an irreducible symmetric space of noncompact type}
We will assume that $X$ has rank at least two, since otherwise
there is nothing to prove.
The strategy of the proof is to use $B$ to produce a subgroup 
$H\subset G$  which has no fixed point in $\tits X$, which 
can be used to tie $B$ closely with $X$.  When $B$ is irreducible,
$H$ is generated using ``restricted'' root groups, and when
$B$ is reducible $H$ is generated by transvections, and
decomposes as a product.

We let $W$ denote the Weyl group of $X$.  Thus  $\tits X$
is a spherical building modelled on a spherical 
Coxeter complex $(S,W)$.  We let $W_B\subset W$ denote the
sub-Coxeter group defining a {\em thick} building structure
on $B$, see \cite[sec. 3.7]{klle}; thus each $W_B$-wall
in $B$ lies in at least $3$ roots (or half-apartments)
of $B$.

\medskip
{\em Case 1.  The subbuilding $B$ is irreducible.}
Our first step is to show that the Moufang property 
restricts to top dimensional
irreducible subbuildings.  Let $a\subset B$ be a  $W_B$-root in 
$B$.  Let $U_a\subset \aut(\tits X)$ 
denote the root group of $a$ (see section \ref{prelimrootgroups}).

\begin{definition}
The {\em restricted root group} of $a$ is defined to be
the subgroup  $U_a^B\subset U_a$ which preserves the subbuilding 
$B\subset \tits X$.
\end{definition}

\begin{lemma}
$U_a^B$ acts 
transitively on the collection of roots in $B$ opposite to
$a$.
\end{lemma}
\proof  Pick two $W_B$-roots $a_1,a_2\subset B$ opposite
$a$.  Since $\tits X$ is Moufang, there is a unique $g\in U_a$
such that $g(a_1)=a_2$.  Let $B'\defeq B\cap g^{-1}(B)$.
Note that $B'\subset B$ is a convex subset 
(see Definition \ref{defconvex})
containing the apartment $a\cup a_1$; therefore by 
\cite[Prop. 3.10.3]{klle}, $B'$ is a top dimensional subbuilding
of $B$.  Let $\si\subset a$
be a $W$-chamber disjoint from the boundary $\D a$, and for $i=1,2$
let $\si_i\subset a_i$ be the chamber  in
$a_i$ opposite $\si$ ; likewise, let  $\pi\subset \si$ be a panel (a 
codimension $1$ face) of $\si$, and for $i=1,2$ let $\pi_i\subset \si_i$ 
be the opposite panel in $a_i$.  Now for each chamber
$\si'\subset B$ incident to $\si$ along $\pi$, for each 
$i=1,2$ there is a unique chamber $\si_i'$ incident to 
$\si_i$ along $\pi_i$, which corresponds to $\si'$
under the correspondence of \cite[Prop. 3.6.4]{klle}; clearly
$g(\si_1)=\si_2$, and hence $g(\si_1')=\si_2'$.  This implies
that $\si_1'\subset B'$.  Now we may argue as in the proof
of \cite[Prop. 3.12.2]{klle} to see that $B'=B$, and  therefore
$g(B)\subset B$; applying the same reasoning to $g^{-1}$
we conclude that $g(B)=B$.  Thus we have shown that $U_a^B$
acts transitively on the roots in $B$ opposite $a$. 
\qed

\medskip
Now pick a $W_B$-wall $\omega\subset B$, and let 
$\tits X(\omega)\subset\tits X$
be the subbuilding consisting of the union of the apartments containing
$\omega$; similarly, let $B(\omega)$ be the subbuilding of $B$
determined by $\omega$.  Thus if $F\subset X$ is a singular 
flat with $\tits F=\omega$,
then the parallel set $\P(F)$ has Tits boundary $\tits X(\omega)$,  
the product
splitting $\P(F)=F\times Y$ induces a join decomposition
$\tits X(\omega)=\omega\circ \;\tits Y$, and $Y\subset X$ is a rank $1$ 
symmetric 
subspace of dimension $>1$.   This join decomposition induces
a join decomposition $B(\omega)=\omega\circ \La$, where
$\La\defeq \tits Y\cap B$.

\begin{lemma}
$\La$ is a compact connected manifold of positive dimension.
\end{lemma} 
\proof
We observe that for each root $a\subset \tits X$
with $\D a=\omega$, the root group $U_a$ acts freely
transitively by homeomorphisms on $\geo Y\setminus \{\xi\}$, 
where $a=\omega\circ \xi$.  Thus if we choose 
$\xi'\in \tits Y\setminus \{\xi\}$ and let $a'\defeq\omega\circ \xi'$,
 then the 
map $\phi:U_a\ra \geo Y\setminus\{\xi\}$
defined by $\phi(g)\defeq g\xi'$
is a continuous bijection between manifolds, and is therefore
a homeomorphism.    Now suppose  $\xi,\xi'\in \La$, so that
$a,a'\subset B$.  The restricted root group $U_a^B\subset U_a$ 
acts simply transitively
on $\La\setminus \{\xi\}$, so 
$\phi$ restricts to a homeomorphism $U_a^B\ra \La\setminus\{\xi\}$.
Thus $U_a^B$ is a closed subgroup of $U_a$, and is therefore a 
manifold, which means that  $\La\setminus\{\xi\}$ is also a manifold.   
Since 
$\xi\in \La$ was chosen arbitrarily, it follows that the 
group generated by the collection of restricted root
groups $\{U_a^B\mid a=\omega\circ\;\xi,\;\;\xi\in \La\}$,
acts transitively on $\La$.  Thus $\La$ is a compact manifold.

Note that $|\La|\geq 3$, since $\La$ is in bijection
with the roots of $B$ containing $\omega$.  Since $U_a$ 
is unipotent, every $g\in U_a^B\setminus\{e\}$ has infinite
order.  This implies that $\La$ is an infinite set; being
a compact manifold, it must have positive dimension.

If $\xi\in \La$ and $a\defeq \omega\circ\xi$, then $U_a^B$
acts transitively on $\La\setminus\{\xi\}$ while preserving
the connected component of $\La$ containing $\xi$.  It follows
that $\La$ is connected.  
\qed

\medskip
Let $H$ be the subgroup of $G$ generated by the 
restricted root groups $U_a^B$, where $a$
ranges over all $W_B$-roots in $B$.  $H$ is a connected subgroup
of the Lie group $G$ since it is generated by
connected subgroups.

Our next objective is to show that $H$ does not fix any point
in $\tits X$.

\begin{lemma}
\label{foldchamber}
Suppose $\si_0\subset\tits X$ is a chamber intersecting
$B$ in a panel $\pi_0$.  Then there is a sequence $h_k\in H$
such that $h_k\si_0$ converges in $\geo X$ to a chamber
$\si_1\subset B$.
\end{lemma}
\proof
We may assume that $\si_0\not\subset B$.  Therefore $\pi_0=\si_0\cap B$
is contained in a unique $W_B$-panel $\pi_1\subset B$.  Hence there
is a $W_B$-chamber $\si_1\subset B$ containing $\pi_1$, and 
$\pi_1$ is contained in a $W_B$-wall $\omega\subset B$.
The chambers $\si_0$ and $\si_1$ lie in roots $a_0\subset\tits X(\omega)$,
and $a_1\subset B(\omega)$ respectively.    Using the notation above, we have
join decompositions $\tits X(\omega)=\omega\circ\;\tits Y$,
$\geo X(\omega)=\omega\circ\;\geo Y$.   For $i=0,1$ let 
$\xi_0,\xi_1\in \geo Y$
be the element such that $a_i=\omega\circ \xi_i$.  The unipotent
root group $U_{a_1}^B$ acts by parabolic homeomorphisms on $\geo Y$ fixing
$\xi_1\in \geo Y$.  If $g\in U_{a_1}^B\setminus\{e\}$, it has
infinite order, and $g^k\xi_0$ converges as $|k|\ra\infty$
to $\xi_1$ with respect to the topology of $\geo Y$.  It follows
that $g^k\si_0$ converges to $\si_1$ as $|k|\ra \infty$.
\qed

\begin{lemma}
\label{collapsetob}
For every chamber $\si_1\subset \tits X$, there is a sequence
$h_k\in H$ such that $h_k\si_1$ converges to a chamber in $B$.
\end{lemma}
\proof
Let $\si_1,\si_2,\ldots,\si_k$ be a gallery in the spherical
building $\tits X$, where $\si_k\subset B$.   By Lemma \ref{foldchamber},
the Lemma holds if $k=1$, so assume $k>1$, pick $1<j\leq k$,
and suppose the lemma holds for $\si_{j-1}$.  Thus there
is a sequence $g_k\in H$ such that $g_k\si_{j-1}$ converges
to a chamber $\tau\subset B$, and after passing to a subsequence,
we may assume that $g_k\si_j$ converges to a chamber $\tau'\subset
\tits X$ meeting $\tau$ along a $W$-panel $\pi$ (at least). 
We are done if $\tau'\subset B$, so we assume $\tau\not\subset B$,
which implies that $\pi$
is contained in a $W_B$-panel $\pi'$. Let $U\subset \geo X$ be an open
subset containing $B$.  Applying Lemma \ref{foldchamber}, 
there is an $h\in H$ such that $h\tau'\subset U$.
Then for large $k$ we have $hg_k\si_j\subset U$.  The open 
set $U$ was arbitrary, so the lemma follows by induction.
\qed

\medskip
The lemma implies that any point in $\tits X$ fixed by $H$
must lie in $B$;  since $H$ acts transitively on
the set of $W_B$-chambers of the irreducible building $B$,
no such fixed point exists.  

Now suppose  $H$ preserved a proper symmetric
subspace $Y\subset X$.
Then $\tits Y\subset \tits X$ would be a proper  $H$-invariant
subbuilding which defined a closed subset of $\geo X$.
Then  Lemma \ref{collapsetob} forces $B\subset\tits Y$,
which contradicts the assumption that $B$ is not contained
in the boundary of a proper symmetric subspace of $X$.
Thus $H$ is a connected subgroup of $G$ which neither
fixes a point in $\tits X$ nor preserves a proper symmetric
subspace of $X$, and so we conclude that $H=G$, see
section \ref{prelimsymmspace}.  Therefore $B=\tits X$.

\medskip
{\em Case 2. The subbuilding $B$ is reducible.}

\begin{lemma} $B$ cannot have a nontrivial
spherical join factor.
\end{lemma}
\proof   Let $S\subset B$ be a maximal spherical
join factor of $B$, and let $F\subset X$ be a
flat with $\tits F=S$.  Then the boundary of the parallel set 
$\P(F)$  contains $B$.  By our assumption we may
conclude that $X=\P(F)$.  However, $X$ is an irreducible symmetric
space of noncompact type, so this is a contradiction.
\qed

 Let $$B=B_1\circ\ldots\circ
B_l$$ 
be the unique join decomposition of $B$
into irreducible nonspherical join factors.
By case 1 above we are done if there is only
one factor, so we assume that $l>1$. 

For each $i$, we let $H_i\subset G$ be the connected
Lie group generated by transvections along geodesics
whose ideal endpoints lie in $B_i$.  Since transvections
along parallel geodesics coincide, and  transvections
along geodesics lying in a single flat commute, it
follows that $H_i$ commutes with $H_j$ when $i\neq j$.

Let $H\defeq H_1\times\ldots\times H_l$.

\begin{lemma}
$H$ does not fix any point in $\tits X$.
\end{lemma}
\proof
Pick a maximal flat $F\subset X$ such that $\tits F\subset B$.
As $H$ contains the full transvection group of $F$, we get
$\fix(H,\tits X)\subset \tits F$.  This means that the fixed
point set of $H$ is contained in the intersection $S$ of the
apartments of $B$; this intersection is empty since $B$
has no spherical join factor.
\qed

 We must therefore have 
$$H=H_1\times\ldots\times H_l=G,$$
see section \ref{prelimsymmspace}.
This contradicts the fact that $G$ is a simple Lie 
group.

\section{Convex sets preserved by Zariski dense groups}

\begin{theorem}
\label{thmconvexzdense}
Let $X$ be a symmetric space of noncompact type
with de Rham decomposition $X=X_1\times\ldots\times X_k$,
let $\pi_i:X\ra X_i$ be the projection map, 
and    $G= \isom_o(X)$ be the associated
connected semi-simple Lie group.  We denote
by $X=Y_1\times Y_{\geq 2}$ the decomposition
of $X$ into (the product of the) rank $1$ and 
the higher rank factors.
Suppose $\Gamma\subset G$ is a Zariski dense subgroup
which preserves a closed convex 
subset $C\subset X$.  Then $C$ is of the form
\begin{equation}
\label{csplits}
C_1\times Y_{\geq 2},
\end{equation}
where $C_1\subset Y_1$ is convex. 
Furthermore,  for each de Rham factor $X_i$ of
$Y_1$, there
is a  $\Ga$-invariant subset $\hat C_i\subset X_i$  such that
\begin{itemize}
\item $\hat C_i$ is the 
closed convex hull of its limit set.

\item $|\geo C_i|=\infty$,

\item $\hat C_1\defeq \prod_i\hat C_i\subset C_1$.

\item  $\geo\hat C_1=\geo C_1$. 
\end{itemize}
\end{theorem}
\proof
By Lemma \ref{lemtitsconvex}, the limit set
$\La(C)=\geo C$ is a (cone
topology) closed convex
 subset containing
the limit set of $\Gamma$.  By Benoist \cite{benoist},
the limit set of $\Gamma$ contains an open neighborhood
(with respect to the topology of $\tits X$) of a pair of antipodal
regular points $\xi,\,\hat\xi\in\tits X$.  Hence
$\tits C$ contains an apartment in $\geo X$.  By
\cite[Prop. 3.10.3]{klle} it follows that $\tits C$ is a 
top dimensional subbuilding of $\tits X$.  

Suppose $\tits C\subset \tits Y$ for some proper
symmetric subspace $Y\subset X$.  For every
apartment $A\subset \tits C$, there is a unique maximal
flat $F\subset X$ with $\tits F=A$, and so
$F\subset Y$; likewise, we have $F\subset gY$ for all
$g\in \Ga$ which implies that $F\subset\cap_{g\in\Ga}gY$.
Since $A$ was chosen arbitrarily, we conclude that
$\cap_{g\in\Ga}gY\subset X$ is a $\Ga$-invariant proper
symmetric subspace, which contradicts the Zariski
density of $\Ga$.

  Theorem \ref{filling} applies, so  the Tits boundary $\tits C$ 
splits as a join
$\tits C=B_1\circ\ldots\circ B_k$, where $B_i=\tits X_i$
when $X_i$ has rank at least two, and $|B_i|=\infty$
for each $i$, by the Zariski density of $\Gamma$.  

Applying  Lemma \ref{lemconvexcone}, it follows that 
$C$ splits as in (\ref{csplits}).

Define $\hat C_i\subset X_i$ to
be the closed convex hull of $B_i$; when $\rank(X)\geq 2$ then
$\hat C_i=X_i$.  Applying  Lemma \ref{lemconvexcone},
it follows that $\hat C_1\defeq \prod_i\hat C_i\subset C_1$.
\qed

\section{Invariant convex subsets 
in symmetric spaces with Euclidean deRham factors}

\begin{theorem}
\label{thmconvexrn}
Let $Y$ be a symmetric space 
of noncompact type 
without Euclidean de Rham
factor,
and suppose $\Ga\subset \isom(Y)\times \isom(\E^n)$
is a subgroup
whose projection $\pi_Y(\Ga)\subset \isom(Y)$
is Zariski dense in the identity component $\isom_o(Y)$.
If $C\subset X\defeq Y\times\E^n$
is a $\Ga$-invariant closed convex set,
then either $C=\pi_Y(C)\times \E^n$
or there is a proper $\Ga$-invariant affine subspace $A\subset \E^n$.
\end{theorem}
\proof
We denote by $\sh:=\pi_Y(C)$ 
the shadow of $C$ in $Y$.
For every point $y\in \sh$ we consider the slice 
$(\{y\}\times\E^n) \cap C=:C_y$. 
Since $C$ is closed,
the boundary at infinity $\tits C_y$
does not depend on $y$ 
and it is a closed convex subset $D$ of the round $(n-1)$-sphere $\tits \E^n$. 
We may assume that it is a proper subset
because otherwise $C=\sh\times\E^n$ 
and we are done. 

If the $C_y$ split off an $\R^k$-factor, 
$1\leq k<n$, 
then $C$ itself splits off an $\R^k$-factor. 
If $E'\subset\E^n$ is the maximal Euclidean factor 
and $\E^n=E'\times E''$ a splitting 
then this splitting is preserved by $\Ga$.
We can therefore reduce to the case that the $C_y$ 
have no Euclidean factor. 

{\em Case 1: The slices $C_y$ are unbounded.}
The set $D\subset\tits\E^n$ has diameter $<\pi$ 
and hence a well-defined center $\zeta$
which must be fixed by $\Ga$.
Let $b_{\zeta}$ denote the Busemann function on $X$ associated to $\zeta$.
For every $\ga\in\Ga$ the difference 
$b_{\zeta}(\ga\,\cdot)-b_{\zeta}$ 
equals a constant $\rho(\ga)$
and the map $\rho:\Ga\to\R$ is a group homomorphism. 

The restriction of $b_{\zeta}$ to $C_y$ is bounded above 
because $\tits C_y$ is contained in the open ball 
$B_{\frac{\pi}{2}}(\zeta)$.
We may therefore 
assign to each $y\in \sh$ the 
bottom height of the slice $C_y$ in the direction $\zeta$
defined as 
$h(y):=\min(-b_{\zeta}|_{C_y})$. 
The function $h:\sh\ra\R$ is convex. 
We consider the asymptotic slope function 
$\slope_h:\tits \sh\to \R\cup\{\infty\}$, 
see section \ref{prelimslopes}.
It is $\Ga$-invariant. 
If the homomorphism $\rho$ is nontrivial
then $\slope_h$ assumes also negative values, 
and by Proposition \ref{propuniqueminimum} it
has a unique minimum. 
This minimum must be fixed by $\Ga$,
a contradiction to the Zariski density of $\pi_Y(\Ga)$ in $\isom(Y)$.
Therefore $\rho$ must be trivial, and 
the level sets of $b_{\zeta}$ 
yield $\Ga$-invariant hyperplanes in $\E^n$.

{\em Case 2: The slices $C_y$ are bounded.}
We pick an ideal point $\zeta\in\tits E^n$. 
As above, 
measuring the height in the direction of $\zeta$,
we can consider 
the convex function $\bot:\sh\to\R$
given by 
$\bot(y):=\min(-b_{\zeta}|_{C_y})$ 
and the concave function $\top:\sh\to\R$
given by 
$\top(y):=\max(-b_{\zeta}|_{C_y})$. 
both functions are continuous because $C$ is closed.

We now use the structure Theorem \ref{thmconvexzdense}
for convex sets invariant under a Zariski dense group. 
It implies that $\tits \sh$ splits as the spherical join
of the boundaries of the higher rank factors 
and of infinite subsets in the boundaries of the rank one factors. 
In particular,
$\geo \sh$ has a well-defined 
and therefore $\pi_Y(\Ga)$-invariant 
convex hull $\ch(\geo \sh)$ in $Y$
which is the product of the higher rank factors of $Y$
with the closed convex hulls of the subsets 
in the boundaries of the rank one factors. 

Lemma \ref{concavecomplete} applied to the higher rank factors 
and 
Lemma \ref{nonnegconcrank1} applied to the rank one factors 
imply that the continuous concave function 
$\top-\bot:\sh\to[0,\infty)$
is constant on $\ch(\geo \sh)$. 
It follows that the restrictions of $\top$ and $\bot$ to $\ch(\geo \sh)$
are affine . 
According to 
Lemmas \ref{affinehigherrank} and \ref{affinerank1}
both functions are constant on $\ch(\geo \sh)$. 

Since the values of $\top(y)$ (or $\bot(y)$) for all directions $\zeta$
determine the slice $C_y$
it follows that the slices $C_y$ equal 
the same compact set $B\subset E^n$ 
for all $y$ in the $\pi_Y(\Ga)$-invariant subset $\ch(\geo \sh)$. 
In particular, 
the action of $\Ga$ on $\E^n$ has bounded orbits 
and therefore a fixed point.
\qed 

\section{The convex cocompact case}

In this section we prove:

\begin{lemma}
\label{lemnofixed}
Let $X=\E^n\times Y$, 
where $Y$ is a symmetric space of noncompact type.
If $\Ga \subset \isom(X)$ is a discrete convex
cocompact group which does not preserve any 
proper symmetric subspace of $X$, then the fixed 
point set of $\Ga$ in   $\tits X$ is contained in the
Tits boundary of the Euclidean factor $\E^n$.
\end{lemma}
\proof
Let $C$ be a $\Gamma$-invariant closed convex set
on which $\Ga$ acts cocompactly.  
Suppose $\Ga$ fixes a point $\xi\in\geo X\setminus \geo\E^n$.
The $\Ga$-action respects the join structure of $\tits X$,
so we may assume without loss of generality that $\xi\in\geo Y$.

Recall that since $\Ga$ fixes $\xi$, the $\Ga$-translates
of the Busemann function $b_\xi$ differ by a constant,
and the map $\Ga\ni g\mapsto g_*(b_\xi)-b_\xi$
defines a homomorphism $\rho:\Ga\to\R$.

Suppose first that the homomorphism $\rho$ is trivial,
i.e.\ $b_{\xi}$ is $\Ga$-invariant.  
Then $b_\xi|_C$ is bounded and attains a minimum.
The minimum set of $b_\xi|_C$ is a convex subset $C_1\subset C$ lying
in a horosphere.  By triangle comparison one concludes that
if $p_1,p_2\in C_1$, then the ideal geodesic triangle
$\ol{\xi p_1}\cup\ol{p_1p_2}\cup\ol{p_2\xi}$ bounds a flat
half-strip.  Thus $C_1$ is contained in the parallel 
set $\P(\ga)$ of a geodesic $\ga\subset \E^n\times Y$ which is 
parallel to the $Y$ factor. 
Since $C_1$ is $\Ga$-invariant 
it follows that $\Ga$ preserves
a proper symmetric subspace of $X$, which is a contradiction.
Therefore $\rho$ is a nontrivial homomorphism and $b_\xi(C)=\R$.

Consider a group element $g\in\Ga$ which 
translates the Busemann function $b_{\xi}$.
We may assume that $b_{\xi}(gx)=b_{\xi}(x)-a$ 
for all $x\in X$
with $a>0$.  As the action is discrete, $\Ga$ acts
on $C$ by semi-simple isometries, and so $g$
is an axial isometry.
  Pick a point $x_0\in C$ and let $\rho:[0,\infty)\ra X$ 
be the unit speed ray starting in $x_0$ and asymptotic to $\xi$.
Then for $x_n=g^nx_0$ holds
$b_{\xi}(x_n)=b_{\xi}(\rho(na))$.
We obtain that 
$d(x_n,\rho(na))\leq n\, d(x_1,\rho(a))$ 
and 
$\angle_{\rho(na)}(x_n,x_0)\geq\frac{\pi}{2}$.
Triangle comparison implies that 
$$
\tan\tangle(\xi_1,\xi)\leq\frac{d(x_1,\rho(a))}{a}>0
$$ 
and thus 
$\tangle(\xi_1,\xi)<\frac{\pi}{2}$.

Since $\tangle(\xi,\geo C)<\frac{\pi}{2}$ 
there is a unique $\eta\in\geo C$
at minimum Tits distance from $\xi$, and so $\eta$ is fixed
by $\Gamma$.  As $\tangle(\eta,\xi)< \frac{\pi}{2}$, it
follows that $\eta$ does not lie in $\tits\E^n$.

We now apply
Lemma \ref{prelimrsplit} to the convex set $C$. 
We obtain that the convex set $C$
contains a $\Ga$-invariant parallel set 
(with respect to $C$)
$Z\defeq \P(\ga)\subset C$,
where $\tits\ga\ni \eta$.  Therefore $\Ga$ preserves
the parallel set of $\ga$ in $X$, which is a contradiction.
\qed

\section{The proof of Theorems \ref{mainzdense} and 
\ref{mainconvexcocompact}}

{\em Proof of Theorem \ref{mainzdense}.}
This  follows immediately from 
Theorem \ref{thmconvexrn} and Theorem \ref{thmconvexzdense}.

{\em Proof of Theorem \ref{mainconvexcocompact}.}
By Lemma \ref{lemnofixed} and the fact that $X$ 
contains no proper $\Ga$-invariant symmetric subspace,
the fixed point set of $\Ga$ is contained in $\E^n$.
Therefore the projection of $\Ga$ to $\isom(Y)$ 
is Zariski dense in $\isom_o(Y)$, since otherwise it would preserve
a proper symmetric subspace $Y'\subset Y$, contradicting
our assumption on $X$.  The theorem then follows from 
Theorem \ref{mainzdense}.
\qed

\bibliographystyle{alpha}
\bibliography{convex}

\begin{thebibliography}{And83}

\bibitem[And83]{anderson}
Michael~T. Anderson.
\newblock The {D}irichlet problem at infinity for manifolds of negative
  curvature.
\newblock {\em J. Differential Geom.}, 18(4):701--721 (1984), 1983.

\bibitem[Bal95]{ball}
Werner Ballmann.
\newblock {\em Lectures on spaces of nonpositive curvature}, volume~25 of {\em
  DMV Seminar}.
\newblock Birkh\"auser Verlag, Basel, 1995.
\newblock With an appendix by Misha Brin.

\bibitem[Ben97]{benoist}
Y.~Benoist.
\newblock Propri\'et\'es asymptotiques des groupes lin\'eaires.
\newblock {\em Geom. Funct. Anal.}, 7(1):1--47, 1997.

\bibitem[BH99]{brha}
Martin~R. Bridson and Andr{\'e} Haefliger.
\newblock {\em Metric spaces of non-positive curvature}, volume 319 of {\em
  Grundlehren der Mathematischen Wissenschaften [Fundamental Principles of
  Mathematical Sciences]}.
\newblock Springer-Verlag, Berlin, 1999.

\bibitem[BL04]{ballyt}
A.~Balser and A.~Lytchak.
\newblock Centers of convex subsets of buildings.
\newblock ArXiv preprint math.MG/0410440, 2004.

\bibitem[BT71]{boreltits}
A.~Borel and J.~Tits.
\newblock \'{E}l\'ements unipotents et sous-groupes paraboliques de groupes
  r\'eductifs. {I}.
\newblock {\em Invent. Math.}, 12:95--104, 1971.

\bibitem[KL97]{klle}
Bruce Kleiner and Bernhard Leeb.
\newblock Rigidity of quasi-isometries for symmetric spaces and {E}uclidean
  buildings.
\newblock {\em Inst. Hautes \'Etudes Sci. Publ. Math.}, (86):115--197, 1997.

\bibitem[Mos55]{mos}
G.~D. Mostow.
\newblock Some new decomposition theorems for semi-simple groups.
\newblock {\em Mem. Amer. Math. Soc.}, 1955(14):31--54, 1955.

\bibitem[Qui04]{quint}
J.-F. Quint.
\newblock Groupes convexes cocompacts en rang sup\'erieur.
\newblock preprint, 2004.

\bibitem[Ron89]{ronan}
Mark Ronan.
\newblock {\em Lectures on buildings}, volume~7 of {\em Perspectives in
  Mathematics}.
\newblock Academic Press Inc., Boston, MA, 1989.

\bibitem[Tit74]{tits}
Jacques Tits.
\newblock {\em Buildings of spherical type and finite {BN}-pairs}.
\newblock Springer-Verlag, Berlin, 1974.
\newblock Lecture Notes in Mathematics, Vol. 386.

\end{thebibliography}

\end{document}